\documentclass[a4paper,12pt,reqno]{amsart}

\usepackage[T1]{fontenc}
\usepackage[utf8]{inputenc}
\usepackage{lmodern}
\usepackage[english]{babel}

\usepackage[%
	left=2.5cm,       
	right=2.5cm,      
	top=3.5cm,        
	bottom=3.5cm,     
	heightrounded,    
	bindingoffset=0mm 
]{geometry}

\usepackage{amsmath}
\usepackage{amssymb}
\usepackage{mathtools}
\usepackage{mathrsfs}
\usepackage[normalem]{ulem}

\numberwithin{equation}{section}

\usepackage{hyperref} 
\usepackage[noabbrev,capitalize]{cleveref}

\usepackage{bookmark}

\usepackage[initials,
]{amsrefs}

\theoremstyle{plain}
	\newtheorem{theorem}{Theorem}[section]
	
	\newtheorem{proposition}[theorem]{Proposition}
	\newtheorem{corollary}[theorem]{Corollary}

\theoremstyle{definition}

	\newtheorem{remark}[theorem]{Remark}
	\newtheorem{open.problem}[theorem]{Open Problem}

\newcommand{\R}{\mathbb{R}}

\newcommand{\Sph}{\mathbb{S}} 
\newcommand*{\transp}[2][-3mu]{\ensuremath{\mskip1mu\prescript{\smash{\mathrm t\mkern#1}}{}{\mathstrut#2}}}	

\newcommand{\eps}{\varepsilon}

\usepackage{tikz}
\usepackage{graphicx}
\usepackage{xcolor}
\usepackage[font=small]{caption}

\newcommand{\closure}[2][3]{%
  {}\mkern#1mu\overline{\mkern-#1mu#2}}
  
\usepackage{enumerate}

\usepackage{url}

\newcommand{\de}{\partial}

\renewcommand{\phi}{\varphi}
\renewcommand{\rho}{\varrho}
\renewcommand{\theta}{\vartheta}
\newcommand{\M}{\mathscr{M}}

\DeclareMathOperator{\supp}{supp}
\DeclareMathOperator{\sgn}{sgn}

\DeclareMathOperator{\SO}{SO}

\DeclareMathOperator{\Lip}{Lip}
\DeclareMathOperator{\loc}{loc}
\DeclareMathOperator{\NL}{NL}

\newcommand{\e}{\mathrm{e}}

\DeclarePairedDelimiter{\set}{\{}{\}}

\mathchardef\ordinarycolon\mathcode`\:
\mathcode`\:=\string"8000
\begingroup \catcode`\:=\active
  \gdef:{\mathrel{\mathop\ordinarycolon}}
\endgroup

\def\Xint#1{\mathchoice
{\XXint\displaystyle\textstyle{#1}}%
{\XXint\textstyle\scriptstyle{#1}}%
{\XXint\scriptstyle\scriptscriptstyle{#1}}%
{\XXint\scriptscriptstyle\scriptscriptstyle{#1}}%
\!\int}
\def\XXint#1#2#3{{\setbox0=\hbox{$#1{#2#3}{\int}$ }
\vcenter{\hbox{$#2#3$ }}\kern-.6\wd0}}

\def\aint{\Xint-}

\newcommand{\Haus}[1]{\mathscr{H}^{#1}} 
\newcommand{\Leb}[1]{\mathscr{L}^{#1}} 

\renewcommand{\div}{\mathrm{div}} 
\newcommand{\redb}{\mathscr{F}} 

\newcommand{\res}{\mathop{\hbox{\vrule height 7pt width .5pt depth 0pt
\vrule height .5pt width 6pt depth 0pt}}\nolimits}

\allowdisplaybreaks

\begin{document}

\title{Failure of the local chain rule for the fractional variation}

\author[G.~E.~Comi]{Giovanni E. Comi}
\address[G.~E.~Comi]{Centro di Ricerca Matematica ``Ennio De Giorgi'', Scuola Normale Superiore, Piazza dei Cavalieri 3, 56126 Pisa (PI), Italy}
\email{giovanni.comi@sns.it}

\author[G.~Stefani]{Giorgio Stefani}
\address[G.~Stefani]{Scuola Internazionale Superiore di Studi Avanzati (SISSA), via Bonomea~265, 34136 Trieste (TS), Italy}
\email{giorgio.stefani.math@gmail.com}

\thanks{\textit{Acknowledgments}. 
The authors thank Daniel Spector for many valuable comments on a preliminary version of the present work and for pointing out the references~\cites{M75,RSS22}.
The authors are members of  the Istituto Nazionale di Alta Matematica (INdAM), Gruppo Nazionale per l'Analisi Matematica, la Probabilità e le loro Applicazioni (GNAMPA).
The first author is partially supported by the INdAM--GNAMPA 2022 Project \textit{Alcuni problemi associati a funzionali integrali: riscoperta strumenti classici e nuovi sviluppi}, codice CUP\_E55F22000270001, and has received funding from the MIUR PRIN 2017 Project ``Gradient Flows, Optimal Transport and Metric Measure Structures''.
The second author is partially supported by the INdAM--GNAMPA 2022 Project \textit{Analisi geometrica in strutture subriemanniane}, codice CUP\_E55F22000270001, and has received funding from the European Research Council (ERC) under the European Union’s Horizon 2020 research and innovation program (grant agreement No.~945655). Part of this work was undertaken while the authors were visiting each other at the University of Pisa and the Scuola Internazionale Superiore di Studi Avanzati (SISSA) in Trieste. They would like to thank these institutions for the support and warm hospitality during the visits.
}

\keywords{Fractional gradient, fractional divergence, fractional variation, fractional Hardy inequality, chain rule}

\subjclass[2020]{Primary 46E35. Secondary 28A12}

\date{\today}

\begin{abstract}
We prove that the local version of the chain rule cannot hold for the fractional variation defined in \cite{CS19}. In the case $n = 1$, we prove a stronger result, exhibiting a function $f \in BV^{\alpha}(\R)$ such that $|f| \notin BV^{\alpha}(\R)$. The failure of the local chain rule is a consequence of some surprising rigidity properties for non-negative functions with bounded fractional variation which, in turn, are derived from a fractional Hardy inequality localized to half-spaces. 
Our approach exploits the results of~\cite{CS21} and the distributional approach of the previous papers~\cites{BCCS20,CS19,CS19-2,CSS21}. 
As a byproduct, we refine the fractional Hardy inequality obtained in \cites{SS18, S20-An} and we prove a fractional version of the closely related Meyers--Ziemer trace inequality.
\end{abstract}

\maketitle


\section{Introduction}

\subsection{The fractional variation}
Let $\alpha\in(0,1)$.
The \emph{fractional $\alpha$-gradient}  of a function $f\in \Lip_c(\R^n)$ is defined as
\begin{equation}
\label{eqi:def_frac_nabla}
\nabla^\alpha f(x)
=
\mu_{n,\alpha}
\int_{\R^n}\frac{(y-x)(f(y)-f(x))}{|y-x|^{n+\alpha+1}}\,dy,
\quad
x\in\R^n,	
\end{equation}
where
\begin{equation*}
\mu_{n,\alpha}
=
2^{\alpha}\, \pi^{- \frac{n}{2}}\, \frac{\Gamma\left ( \frac{n + \alpha + 1}{2} \right )}{\Gamma\left ( \frac{1 - \alpha}{2} \right )}
\end{equation*}
is a renormalizing constant controlling the behavior of $\nabla^\alpha$ as $\alpha\to1^-$.
A simple computation (see~\cite{CS19}*{Proposition~2.2} for instance) shows that one can equivalently write $\nabla^\alpha f=\nabla I_{1-\alpha}f$ whenever $f\in C^\infty_c(\R^n)$ (and even for less regular functions, see~\cite{CS19}*{Lemma~3.28(i)} for a more precise statement), where
\begin{equation*}
I_{s} f(x) 
= 
2^{-s} \pi^{- \frac{n}{2}} \frac{\Gamma\left(\frac{n-s}2\right)}{\Gamma\left(\frac s2\right)}
\int_{\R^{n}} \frac{f(y)}{|x - y|^{n - s}} \, dy, 
\quad
x\in\R^n,
\end{equation*}
is the \emph{Riesz potential} of order $s\in(0,n)$.

The literature around the operator $\nabla^\alpha$ has been quickly growing  in the recent years in various research directions.
On the one side, we refer the reader 
to~\cites{lo2021class, rodrigues2019nonlocal,rodrigues2019nonlocal-corr,SS15,SS18,SSS15,SSS18} for well-posedness results concerning solutions of PDEs and minimizers of functionals involving this fractional operator, and to~\cites{BCM20,BCM21,KS21} for the the study of polyconvexity and quasiconvexity in connection with the present fractional setting.
On the other side, the properties of $\nabla^\alpha$ led to the discovery of new (optimal) embedding inequalities~\cites{SSS17,S19,S20-An} and the development of a distributional and asymptotic analysis in this fractional framework~\cites{BCCS20,CS19,CS19-2,Sil19,CSS21,CS21}.
For a general panoramic on the fractional framework, the reader may consult the survey~\cite{S20-New} and the monograph~\cite{P16}.

At least for sufficiently smooth functions, the operator $\nabla^\alpha$ obeys the following natural fractional integration-by-parts formula
\begin{equation}
\label{eqi:ibp}
\int_{\R^n}f\,\div^\alpha\phi \,dx=-\int_{\R^n}\phi\cdot\nabla^\alpha f\,dx,
\end{equation}
where 
\begin{equation*}
\div^\alpha \phi(x)
=
\mu_{n,\alpha}
\int_{\R^n}\frac{(y-x)\cdot(\phi(y)-\phi(x))}{|y-x|^{n+\alpha+1}}\,dy,
\quad
x\in\R^n,	
\end{equation*}
is the \emph{fractional $\alpha$-divergence} of the vector field $\phi\in\Lip_c(\R^n;\R^n)$. 

Equality~\eqref{eqi:ibp} is the fundamental basis of the distributional theory in the present fractional setting developed in the previous papers~\cites{BCCS20,CS19,CS19-2,CSS21,CS21}.
In more precise terms, by imitating the classical definition of $BV$ functions, for a given exponent $p\in[1,+\infty]$, we define the (total) fractional variation of a function $f\in L^p(\R^n)$ as
\begin{equation}
\label{eqi:def_frac_var_p}
|D^\alpha f|(\R^n)
=
\sup\set*{\int_{\R^n}f\,\div^\alpha\phi\,dx : \phi\in C^\infty_c(\R^n;\R^n),\ \|\phi\|_{L^\infty(\R^n;\,\R^n)}\le1}.
\end{equation}	
The above definition naturally gives rise to the linear space of $L^p$ functions with bounded \emph{fractional $\alpha$-variation}
\begin{equation*}
BV^{\alpha,p}(\R^n)=
\set*{f\in L^p(\R^n) : |D^\alpha f|(\R^n)<+\infty}
\end{equation*}
that can be endowed with the norm
\begin{equation*}
\|f\|_{BV^{\alpha,p}(\R^n)}
=
\|f\|_{L^p(\R^n)}+|D^\alpha f|(\R^n),
\quad
f\in BV^{\alpha,p}(\R^n).
\end{equation*}
The resulting normed space is Banach and, moreover, one easily checks that $f\in L^p(\R^n)$ belongs to $BV^{\alpha,p}(\R^n)$ if and only if there exists a finite vector-valued Radon measure $D^\alpha f\in\M(\R^n;\R^n)$, the \emph{fractional $\alpha$-variation measure} of $f$, such that 
\begin{equation*}
\int_{\R^n}f\,\div^\alpha\phi\,dx=-\int_{\R^n}\phi\cdot\,dD^\alpha f
\end{equation*}
for all $\phi\in \Lip_c(\R^n;\R^n)$, see~\cite{CSS21}*{Theorem~3}.

In a very similar way, one can define the \emph{distributional fractional Sobolev space}
\begin{equation*}
S^{\alpha,p}(\R^n)
=
\set*{f\in L^p(\R^n): \nabla^\alpha f\in L^p(\R^n;\R^n)}
\end{equation*}
where $\nabla^\alpha f$ stands for the \emph{distributional fractional $\alpha$-gradient}, see~\cite{CS19}*{Definition~3.9}. 
As proved in~\cite{BCCS20}*{Corollary~1} and in~\cite{KS21}*{Theorem~2.7}, $S^{\alpha,p}(\R^n)=L^{\alpha,p}(\R^n)$ whenever $p\in(1,+\infty)$, where $L^{\alpha,p}(\R^n)$ stands for the \emph{Bessel potential space}.
We refer the reader to~\cite{BCCS20}*{Section~2.1} and to the references therein for an agile account on Bessel potential spaces, and to the discussion in~\cite{CS19}*{Section~3.9} for the relations between $L^{\alpha,p}(\R^n)$ and the \emph{Ga\-gli\-ar\-do--So\-bo\-lev--Slo\-bo\-de\-ckij fractional space} $W^{\alpha,p}(\R^n)$.

The study of the space $BV^{\alpha}(\R^n)=BV^{\alpha,1}(\R^n)$ in the \emph{geometric} regime $p=1$ was  initiated in~\cite{CS19}, also in connection with the naturally associated notion of \emph{fractional Caccioppoli perimeter} (see~\cite{CS19}*{Definition~4.1}), and then further investigated in the subsequent works~\cites{CS19-2,BCCS20}.
The fractional variation of an $L^p$ function for an arbitrary exponent $p\in[1,+\infty]$ has been explored in~\cites{CS19-2,CSS21,CS21}.

Throughout this paper, with a slight abuse of notation (that, however, can be rigorously justified thanks to the analysis done in the previous works~\cites{CS19,CS19-2,BCCS20,CS21,CSS21}), in the integer case $\alpha=1$ we let 
\begin{equation*}
BV^{1,p}(\R^n)
=
\set*{f\in L^p(\R^n) : Df \in\M(\R^n;\R^n)}
\end{equation*}  
be the space of $L^p$ functions, $p\in[1,+\infty]$, with bounded variation.

\subsection{Hardy inequality and chain rule}

Due to the central role played by the classical Hardy inequality in the theory of integer as well as of fractional Sobolev spaces, see~\cite{M18} for an account, in~\cites{SS18}, Shieh and Spector  investigated the validity of the natural analogue of the Hardy inequality in the present fractional setting.
In~\cite{SS18}*{Theorem~1.2}, they proved the validity of the following inequality
\begin{equation}
\label{eqi:hardy_SS}
c_{n,\alpha}\,\int_{\R^n}\frac{|f(x)|}{|x|^\alpha}\,dx
\le 
\int_{\R^n}|\nabla^\alpha |f||\,dx
\end{equation}
for all measurable functions $f$ such that $\nabla^\alpha|f|=\nabla I_{1-\alpha}|f|\in L^1(\R^n;\R^n)$, where $c_{n,\alpha}>0$ is a constant depending on $\alpha\in(0,1)$ and $n\ge2$ only.
Actually, the validity of~\eqref{eqi:hardy_SS} for $n=1$ is not explicitly shown in~\cite{SS18}, but one can still recover it via an \textit{ad hoc} modification of their argument.

Motivated by~\eqref{eqi:hardy_SS}, the authors in~\cite{SS18} asked if it is possible to remove the modulus in the right-hand side of~\eqref{eqi:hardy_SS}, that is, more generally, if the following chain rule for the fractional gradient
\begin{equation}
\label{eqi:chain_SS}
\int_{\R^n}|\nabla^\alpha |f||\,dx
\le 
c_{n,\alpha}
\int_{\R^n}|\nabla^\alpha f|\,dx
\end{equation}
holds whenever $f$ is measurable with $\nabla I_{1-\alpha} f\in L^1(\R^n;\R^n)$, where $c_{n,\alpha}>0$ is a constant depending on $\alpha$ and $n$ only, see~\cite{SS18}*{Open Problem~1.4}.

Later, Spector proved the validity of the fractional Hardy inequality
\begin{equation}
\label{eqi:hardy_spector}
c_{n,\alpha}\,\int_{\R^n}\frac{|f(x)|}{|x|^\alpha}\,dx
\le 
\int_{\R^n}|\nabla^\alpha f|\,dx
\end{equation}  
for $n\ge2$, whenever $f\in L^p(\R^n)$ with $p\in\left[1,\frac n{1-\alpha}\right)$ and $\nabla^\alpha f=\nabla I_{1-\alpha} f\in L^1(\R^n;\R^n)$, see~\cite{S20-An}*{Theorem~1.4}.
The approach used in~\cite{S20-An} completely bypasses the validity of~\eqref{eqi:chain_SS} and instead relies on an optimal embedding in Lorentz spaces for the Riesz potential, see~\cite{S20-An}*{Theorem~1.1}.

The relation between the Hardy inequality in~\eqref{eqi:hardy_SS}, as well as the one in~\eqref{eqi:hardy_spector}, with the one valid in the usual fractional Sobolev space $W^{\alpha,1}(\R^n)$ easily follows from the elementary inequality 
\begin{equation*}
\|\nabla^\alpha f\|_{L^1(\R^n;\R^n)}
\le 
\mu_{n,\alpha}\,[f]_{W^{\alpha,1}(\R^n)}
\end{equation*}
naturally available for all functions $f\in W^{\alpha,1}(\R^n)$, see~\cite{CS19}*{Section~1.1}. 
Similar considerations can be done for the Hardy inequalities in the integrability regime $p\in(1,+\infty)$, see the introductions of~\cites{SS15,SS18}.

Up to our knowledge, the validity of  a chain rule for the fractional gradient $\nabla^\alpha$ like~\eqref{eqi:chain_SS} is still an open problem.
Somehow complementing the validity of~\eqref{eqi:hardy_spector} for $n\ge2$, in the present work we disprove the validity of~\eqref{eqi:chain_SS} in the case $n=1$.
More precisely, we prove the following result.     

\begin{theorem}[Failure of the chain rule for $n=1$]
\label{res:chain_sberla_debole}
Let $\alpha\in(0,1)$.
The function 
\begin{equation*}
f_{\alpha}(x) = \mu_{1, - \alpha} \left ( |x|^{\alpha - 1} \sgn{x} - |x - 1|^{\alpha - 1} \sgn(x - 1) \right ),
\quad
x\in\R\setminus\set*{0,1},
\end{equation*}
is such that $f_\alpha\in BV^\alpha(\R)$ but $|f_\alpha|\notin BV^\alpha(\R)$.
\end{theorem}

The proof of \cref{res:chain_sberla_debole} works by contradiction. Precisely, if $D^\alpha |f_\alpha| \in \M(\R)$, then a generalized version of inequality~\eqref{eqi:hardy_SS} for $n = 1$ would hold (see \cref{thm:Hardy} below). 
Thus we would get $f_{\alpha} \in L^1(\R; |x|^{-\alpha} \Leb{1})$, which is clearly false.
Actually, inequality~\eqref{eqi:hardy_SS} cannot be directly applied to the function $f_\alpha$ in \cref{res:chain_sberla_debole}, since $D^\alpha f_\alpha=\delta_0-\delta_1\notin L^1(\R)$, see~\cite{CS19}*{Theorem~3.26}.
However, one can exploit the regularization properties of $BV^{\alpha,p}$ functions~\cite{CSS21}*{Theorem~4} to suitably extend~\eqref{eqi:hardy_SS}, as well as~\eqref{eqi:hardy_spector}, to this more general framework.

\begin{theorem}[Hardy inequality in $BV^{\alpha,p}(\R^n)$ for $p\in{\left[1,\frac n{1-\alpha}\right)}$] \label{thm:Hardy}
Let $\alpha \in (0, 1)$ and $p \in \left [1, \frac{n}{1 - \alpha} \right )$. 
If $f \in BV^{\alpha, p}(\R^n)$, with $f\ge0$ if $n=1$, then
\begin{equation*}
c_{n,\alpha}
\int_{\R^n} \frac{|f(x)|}{|x - x_0|^{\alpha}} \, dx \le |D^\alpha f|(\R^n)
\end{equation*}
for all $x_0 \in \R^n$, where $c_{n,\alpha}>0$ is a constant depending on~$n$ and~$\alpha$ only. In particular, if $n = 1$, the optimal constant is $c_{1, \alpha} = \frac{2 \mu_{1, \alpha}}{\alpha}$.
\end{theorem}

\subsection{Local chain rule}

We do not know if a counterexample to the chain rule~\eqref{eqi:chain_SS} like the one in \cref{res:chain_sberla_debole} can be provided also for $n\ge2$.

The current lack of a counterexample to~\eqref{eqi:chain_SS} may suggest that a stronger version of the chain rule could be valid for the fractional variation for $n\ge2$, in analogy with the chain rule available for $BV$ functions. 
More precisely, for a given $\Phi\in\Lip(\R)$ such that $\Phi(0)=0$, one may wonder if the \emph{local} chain rule 
\begin{equation}
\label{eqi:chain_dura}
|D^\alpha\Phi(f)|\le 
C(\Phi)\,|D^\alpha f|
\quad
\text{in}\ \M(\R^n)
\end{equation}
holds for all $f\in BV^{\alpha,p}(\R^n)$ with $n\ge2$, where $C(\Phi)>0$ is a constant depending on the chosen function $\Phi$ only.
In the present work, we disprove the validity of~\eqref{eqi:chain_dura} for all $n\ge2$ and, actually, we prove the following stronger result.

\begin{theorem}[Failure of the local chain rule]
\label{res:chain_sberla_forte}
Let $\alpha\in(0,1)$ and $p\in\left[1,\frac n{n-\alpha}\right)$.
Let $\Phi\in\Lip(\R)$ be such that $\Phi(0)=0$ and $\Phi\ge0$. 
If $\Phi(f)\in BV^{\alpha,p}(\R^n)$ with
\begin{equation*}
\supp|D^\alpha\Phi(f)|\subset\supp|D^\alpha f|
\end{equation*}
for all $f\in BV^{\alpha,p}(\R^n)$, then $\Phi\equiv0$.
\end{theorem}

In particular, if we consider $\Phi(t) = |t|$ for $t\in\R$, an immediate consequence of \cref{res:chain_sberla_forte} is that, for all $\alpha\in(0,1)$ and $p\in\left[1,\frac n{n-\alpha}\right)$, there exists a function $f \in BV^{\alpha,p}(\R^n)$ such that $\supp|D^\alpha|f||$ is not contained in $\supp|D^\alpha f|$.

The validity of \cref{res:chain_sberla_forte} is a simple consequence again of the analysis made in~\cite{CS19} and of a new surprising rigidity property of non-negative $BV^{\alpha,p}$ functions with $p\in\left[1,\frac n{n-\alpha}\right)$, see \cref{res:null_function} below.
Here and in the rest of the paper, given $\nu \in \Sph^{n - 1}$ and $x_0 \in \R^n$, we let 
\begin{align*}
H_{\nu}^{+}(x_0)  = \{ y \in \R^{n} : (y - x_0) \cdot \nu > 0 \}
\end{align*}
and
\begin{align*}
H_{\nu}(x_0)  = \{ y \in \R^{n} : (y - x_0) \cdot \nu = 0 \}.
\end{align*}
In the case $x_0 = 0$, we simply write $H_{\nu}^{+} = H_{\nu}^{+}(0)$ and $H_{\nu} = H_{\nu}(0)$.
Moreover, for $\alpha\in(0,1)$ and $p\in[1,+\infty]$, we let
\begin{equation*}
BV^{\alpha,p}_+(\R^n)
=
\set*{f\in BV^{\alpha,p}(\R^n) : f\ge0}.
\end{equation*}

\begin{theorem}[Rigidity property in $BV^{\alpha,p}_+(\R^n)$ for $p\in{\left[1,\frac n{n-\alpha}\right)}$]
\label{res:null_function}
Let $\alpha\in(0,1)$, $p\in\left[1,\frac n{n-\alpha}\right)$ and $f\in BV^{\alpha,p}_+(\R^n)$. 
If
either
\begin{equation}
\label{eq:rigidity_support}
\text{$\supp|D^\alpha f|$ is bounded},
\end{equation}
or 
\begin{equation}
\label{eq:rigidity_half-space_+}
|D^\alpha f|\big(\closure[0]{H_\nu^+(x_0)}\big) = 0\
\text{for some $x_0\in\R^n$, $\nu\in\Sph^{n-1}$},
\end{equation}
or
\begin{equation}
\label{eq:rigidity_half-space}
f\in L^\infty(\R^n)\
\text{and}\
D^\alpha f(H^+_\nu(x_0))=0\
\text{for some $x_0\in\R^n$, $\nu\in\Sph^{n-1}$},
\end{equation}
then $f\equiv0$.
\end{theorem}

The rigidity property given by \cref{res:null_function} strongly underlines the difference between the non-local operator~$\nabla^\alpha$ and its local integer counterpart~$\nabla$.
Indeed, it is easily seen that $BV^{1,p}$ functions do not possess such a rigidity property for any given $p\in[1,+\infty]$, due to the locality of the classical variation measure (for instance, one may consider the characteristic function of the unit ball).

In addition, we recall that, despite of the non-local nature of the fractional gradient, there exist functions $f \in BV^{\alpha, p}(\R^n)$, for $p\in\left[1,\frac n{n-\alpha}\right)$, such that $|D^\alpha f|$ is a finite Radon measure with compact support, see the function defined in \cref{res:chain_sberla_debole} for $n = 1$, and \cite{CS19}*{Lemma 3.28} as well as~\cite{CSS21}*{Proposition~4} for the general case. Hence, \cref{res:null_function} immediately tells us that such functions cannot have constant sign.
Conversely, as observed in~\cite{KS21}*{Section~2.2} in the case $n = 1$, given any non-zero function $f\in C^\infty_c(\R^n)$ with $f\ge0$ and $\supp f\subset (-L, L)^n$ for some $L > 0$, for each $j \in \{1, \dots, n\}$ we have $$\nabla^\alpha_j f(x) = {\rm e}_j \cdot \nabla^{\alpha} f(x) \ne0 \text{ at each } x \in \R^n \text{ with } |x_j| \ge L,$$ where ${\rm e}_j$ is the $j$-th vector of the standard coordinate basis of $\R^n$.

We end this section by stating a simple consequence of \cref{thm:Hardy}. To this purpose, we define
$$LSC_{b}(\R^n) = \set*{ f : \R^n \to \R : f \text{ lower semicontinuous and bounded}}.$$

\begin{corollary}
Let $\alpha\in(0,1)$ and $p \in \left [1, \frac{n}{1 - \alpha} \right )$. The operator $$I_{n - \alpha} : BV^{\alpha, p}_+(\R^n) \to LSC_{b}(\R^n)$$ is continuous. In addition, if $n \ge 2$, then $I_{n - \alpha} : BV^{\alpha, p}(\R^n) \to L^{\infty}(\R^n)$ is continuous.
\end{corollary}

\subsection{Integration-by-parts formulas}

The rigidity property of non-negative $BV^{\alpha,p}$ functions stated in 
\cref{res:null_function} is, in turn, a consequence of a fractional Gauss--Green formula on half-spaces, see~\cref{thm:IBP_halfspace} below, which can be regarded as a `vectorial'  Hardy-type equality for the fractional variation. 
Here and in the rest of the paper, for $\alpha\in(0,1)$ and $p,q\in[1,+\infty]$, we let
\begin{equation*}
B^\alpha_{p,q}(\R^n)
=
\set*{
u\in L^p(\R^n)
:
[u]_{B^\alpha_{p,q}(\R^n)}
<
+\infty
}
\end{equation*}
be the space of \emph{Besov functions} on $\R^n$, see~\cite{L17}*{Chapter~17} for its precise definition and main properties, where 
\begin{equation*}
[u]_{B^\alpha_{p,q}(\R^n)}
=
\begin{cases}
\left(
\displaystyle\int_{\R^n}
\frac{\|u(\cdot+h)-u\|_{L^p(\R^n)}^q}{|h|^{n+q\alpha}}
\,dh
\right)^{\frac1q}
&
\text{if}\
q\in[1,+\infty),
\\[6mm]
\sup\limits_{h\in\R^n\setminus\set*{0}}
\dfrac{\|u(\cdot+h)-u\|_{L^p(\R^n)}}{|h|^{\alpha}}
&
\text{if}\
q=+\infty.
\end{cases}
\end{equation*}

\begin{theorem}[Fractional Gauss--Green formula on half-spaces] \label{thm:IBP_halfspace}
Let $\alpha \in (0, 1)$, $p\in\left[1,\frac n{n-\alpha}\right)$ and $q \in \left ( \frac{n}{\alpha}, + \infty \right ]$ be such that $\frac{1}{p} + \frac{1}{q} = 1$.
If $f \in BV^{\alpha,p}(\R^n) \cap L^{\infty}(\R^n)$, then
\begin{equation}
\label{eq:IBP_halfspace}
\frac{\mu_{1,\alpha}}{\alpha}
\lim_{R\to+\infty}
\int_{\R^n} \eta_R(x)\,\frac{f(x)\,\nu}{|(x - x_0) \cdot \nu|^{\alpha}} \, dx
=
- 
\eta(0)
\,
D^{\alpha} f(H_{\nu}^+(x_0))
\end{equation}
whenever $\nu \in \Sph^{n-1}$ and $x_0 \in \R^n$, where $\eta_R(x)=\eta\left(\frac xR\right)$ for $x\in\R^n$ and $R>0$, for some fixed $\eta\in B^\alpha_{q,1}(\R^n)$ with compact support.
In particular, if either $\supp f$ is bounded or $f$ has constant sign, then 
\begin{equation}
\label{eq:IBP_halfspace_special}
\frac{\mu_{1,\alpha}}{\alpha}
\int_{\R^n} \frac{f(x)\,\nu}{|(x - x_0) \cdot \nu|^{\alpha}} \, dx
=
- 
D^{\alpha} f(H_{\nu}^+(x_0)).
\end{equation}
\end{theorem}

We let the reader note that the requirement that $\eta\in B^\alpha_{q,1}(\R^n)$ naturally comes from the general integration-by-parts formula obtained in~\cite{CS21}*{Theorem~1.1}, see~\eqref{eq:leibniz_X} below for a more detailed account.
 
Actually, \cref{thm:IBP_halfspace} is a particular case of the following result, which can be seen as an extension of the integration-by-parts formula~\eqref{eqi:ibp} in the spirit of the fractional Gauss--Green formulas established in~\cite{CS21}*{Section 3.3}.
Here and in the following, we let 
\begin{equation}
\label{eqi:def_precise_repres}
f^\star(x)
=
\begin{cases}
\displaystyle
\lim_{r\to0^+}
\aint_{B_r(x)} f(y)\,dy
&
\text{if the limit exists},
\\[3mm]
0
&
\text{otherwise},
\end{cases}
\end{equation}
be the \emph{precise representative} of $f\in L^1_{\loc}(\R^n)$.

\begin{theorem}[Limit integration-by-parts formula]
\label{res:lim_GG_formula}
Let $\alpha\in(0,1)$ and let  $p\in\left[1,\frac n{n-\alpha}\right)$ and $q\in\left(\frac n\alpha,+\infty\right]$ be such that $\frac1p+\frac1q=1$. 
If $f\in BV^{\alpha,p}(\R^n)\cap L^\infty(\R^n)$ and $g\in W^{\alpha,1}_{\loc}(\R^n)\cap L^\infty(\R^n)$ with $\nabla^\alpha g\in L^1_{\loc}(\R^n;\R^n)$, then
\begin{equation*}
\lim_{R\to+\infty}
\int_{\R^n} \eta_R\, f\,\nabla^\alpha g\,dx
=
-
\eta(0)
\int_{\R^n}g^\star\,dD^\alpha f,
\end{equation*}
where $\eta_R$ is as in \cref{thm:IBP_halfspace} and the limit in \eqref{eqi:def_precise_repres} defining $g^{\star}(x)$ exists for $|D^\alpha f|$-a.e. $x \in \R^n$.
\end{theorem}

In order to apply \cref{res:lim_GG_formula} to get \cref{thm:IBP_halfspace}, one then just need to explicitly compute the fractional gradient of the characteristic function of a half-space. 

\begin{proposition}[$\nabla^\alpha$ of a half-space] \label{prop:halfspace_nabla_alpha}
Let $\alpha \in (0, 1)$, $\nu \in \Sph^{n - 1}$ and $x_0 \in \R^n$. 
We have
\begin{equation} \label{eq:halfspace_nabla_alpha}
\nabla^{\alpha} \chi_{H_{\nu}^{+}(x_0)}(x) 
=
\frac{\mu_{1,\alpha}}{\alpha}
\,
\frac{\nu}{|(x - x_0) \cdot \nu|^{\alpha}} 
\end{equation}
for $x\in\R^n \setminus H_{\nu}(x_0)$.
\end{proposition}

It is worth noticing that \cref{thm:IBP_halfspace} immediately implies the following version of the fractional Hardy inequality for non-negative $BV^{\alpha,p}$ functions in the regime $p\in\left[1,\frac n{n-\alpha}\right)$, where the right-hand side does not involve the knowledge of the  fractional variation on the whole space, but just on a specific half-space.

\begin{corollary}[Fractional Hardy inequality in $BV^{\alpha,p}_{+}(\R^n)$ for $p\in\left[1,\frac n{n-\alpha}\right)$]
\label{res:hardy_mezza}
Let $\alpha \in (0, 1)$ and $p \in \left [1, \frac{n}{n - \alpha} \right )$. 
If $f \in BV^{\alpha, p}_+(\R^n)$, then
\begin{equation} \label{eq:hardy_mezza}
\frac{\mu_{1,\alpha}}{\alpha}
\int_{\R^n} \frac{f(x)}{|(x - x_0)\cdot\nu|^{\alpha}} \, dx \le |D^\alpha f|\big(\closure[0]{H_\nu^+(x_0)}\big)
\end{equation}
for all $x_0 \in \R^n$ and $\nu\in\Sph^{n-1}$.
\end{corollary}

As the reader may notice, \cref{thm:IBP_halfspace} allows to prove \cref{res:null_function} under the assumption~\eqref{eq:rigidity_half-space}.
To deal with the assumption~\eqref{eq:rigidity_support}, one needs to perform a further integration with respect to the direction $\nu\in\Sph^{n-1}$ and obtain the following fractional weighted inequality of Hardy-type.
Again, we underline that the fractional variation appearing in the right-hand side is not computed on the whole space, but just on the complement of a particular ball.

\begin{corollary}[Weighted fractional Hardy-type inequality]
\label{res:weighted_hardy}
Let $\alpha\in(0,1)$ and $p\in\left[1,\frac n{n-\alpha}\right)$. 
If $f\in BV^{\alpha,p}_+(\R^n)$, then 
\begin{equation*}
\int_{\R^n}
f(x)\,w_{n,\alpha}(|x-x_0|,r)\,dx
\le 
|D^\alpha f|\big(\R^n\setminus B_r(x_0)\big)
\end{equation*}
for all $x_0\in\R^n$ and $r>0$, where
\begin{equation*}
w_{n,\alpha}(t,r)
=
\begin{cases}
\displaystyle
\frac{(n-1)\,\omega_{n-1}}{n\,\omega_n}\,\frac{\mu_{1,\alpha}}\alpha
\int_{-1}^{1}\frac{(1-s^2)^{\frac{n-3}2}}{|st-r|^\alpha}\,ds
& 
\text{for}\ n\ge2,
\\[4mm]
\displaystyle
\frac{\mu_{1,\alpha}}{2\alpha}\left(\frac1{|t-r|^\alpha}+\frac1{|t+r|^\alpha}\right)
& 
\text{for}\ n=1.
\end{cases}
\end{equation*}
\end{corollary}

In the particular geometric case $f=\chi_E$ for some measurable set $E\subset\R^n$, the above results read as follows (recall that, by~\cite{CS19}*{Corollary~5.4}, if $\chi_E\in BV^\alpha(\R^n)$, then we have $|D^\alpha\chi_E|\ll\Haus{n-\alpha}$).
Here and in the following, $\redb^\alpha E$ denotes the \emph{fractional reduced boundary} in the sense of De Giorgi, see~\cite{CS19}*{Definition~4.7}.  

\begin{corollary}[Geometric case]
Let $\alpha\in(0,1)$.
If $\chi_E\in BV^\alpha(\R^n)$, then
\begin{align*}
\frac{\mu_{1,\alpha}}{\alpha}
\int_{E} \frac{\nu}{|(x - x_0) \cdot \nu|^{\alpha}} \, dx
&=
- 
D^{\alpha} \chi_E(H_{\nu}^+(x_0)),
\\[3mm]
\int_E |\nabla^{\alpha} \chi_{H_{\nu}^+(x_0)}| \, dx 
&\le 
|D^\alpha \chi_E|(H_{\nu}^+(x_0)),
\\[3mm]
\int_{E}
w_{n,\alpha}(|x-x_0|,r)\,dx
&\le 
|D^\alpha \chi_E|\big(\R^n\setminus B_r(x_0)\big),
\end{align*}
for $x_0\in\R^n$, $\nu\in\Sph^{n-1}$ and $r>0$, where $w_{n,\alpha}$ is as in \cref{res:weighted_hardy}.
Moreover, if either $\supp|D^\alpha\chi_E|$ is bounded or $D^\alpha \chi_E(H^+_\nu(x_0))=0$ for some $x_0\in\R^n$, $\nu\in\Sph^{n-1}$, then $|E|=0$.
In particular, if $|E|>0$, then $\redb^\alpha E$ must be unbounded and must intersect all half-spaces.
\end{corollary}

\subsection{Fractional Meyers--Ziemer trace inequalities}

As discussed in~\cite{S20-New}, the Hardy inequality in~\eqref{eqi:hardy_spector} can be also seen as a particular consequence of known interpolation inequalities in Lorentz spaces.
Precisely, one recognizes that 
\begin{equation*}
\frac1{|\cdot|^{\alpha}}
\in
L^{\frac n\alpha,\infty}(\R^n),
\end{equation*} 
so that~\eqref{eqi:hardy_spector} follows by combining the H\"older inequality
\begin{equation*}
\int_{\R^n}\frac{|f(x)|}{|x|^\alpha}\,dx
\le 
\|f\|_{L^{\frac n{n-\alpha},1}(\R^n)}
\,
\left\|\frac1{|\cdot|^{\alpha}}\right\|_{L^{\frac n{\alpha},\infty}(\R^n)}
\end{equation*}
with the bound
\begin{equation*}
\|f\|_{L^{\frac n{n-\alpha},1}(\R^n)}
\le 
c_{n,\alpha}\,|D^\alpha f|(\R^n)
\end{equation*}
valid for $n\ge2$, which, in turn, is a consequence of~\cite{S20-An}*{Theorem~1.1}.

In the classical integer case, an even more general approach is possible. 
Indeed, if $f\in BV(\R^n)$ and $\mu\in\M_{\loc}^+(\R^n)$ is a non-negative locally finite measure, then
\begin{equation}
\label{eqi:MZ_trace}
\int_{\R^n}|f^\star|\,d\mu
\le c_n\,\|\mu\|_{n-1}\,|Df|(\R^n),
\end{equation} 
for a dimensional constant $c_n>0$, where $f^\star$ is as in~\eqref{eqi:def_precise_repres} and 
\begin{equation*}
\|\mu\|_s=\sup_{x\in\R^n,\ r>0}\frac{\mu(B_r(x))}{r^s}
\end{equation*}
whenever $s\in[0,n]$.
The inequality in~\eqref{eqi:MZ_trace} can be found in~\cite{MZ77}*{Theorem~4.7} and is nowadays called the \emph{Meyers--Ziemer trace inequality}.
We also refer the reader to the recent work~\cite{RSS22} for an interesting historical panoramic around the inequality~\eqref{eqi:MZ_trace}. 
In particular, the authors of~\cite{RSS22} note that V.~G.~Maz\cprime ya proved such an inequality in~\cite{M75}, a few years before the aforementioned~\cite{MZ77}.

Inequality~\eqref{eqi:MZ_trace} plays a central role in the classical $BV$ framework, since it can be considered as the mother inequality of several embedding inequalities, like the Hardy inequality
\begin{equation}
\label{eqi:hardy_bv}
\int_{\R^n}\frac{|f(x)|}{|x|}\,dx
\le c_n\,|Df|(\R^n),
\end{equation}
the Gagliardo--Nirenberg--Sobolev inequality
\begin{equation}
\label{eqi:GNS_bv}
\|f\|_{L^{\frac n{n-1}}(\R^n)}
\le c_n\,|Df|(\R^n),
\end{equation}
and its refinement, the Alvino inequality
\begin{equation}
\label{eqi:alvino_bv}
\|f\|_{L^{\frac n{n-1},1}(\R^n)}
\le c_n\,|Df|(\R^n).
\end{equation}
For a more detailed discussion, we refer the reader to~\cite{S19}*{Section~1} and~\cite{S20-New}*{Section~6}.
Indeed, as soon as $g\in L^{n,\infty}(\R^n)$, one immediately recognizes that the measure
\begin{equation*}
\mu(A)=\int_A g(x)\,dx, 
\quad
A\subset\R^n,
\end{equation*} 
satisfies 
\begin{equation*}
\||\mu|\|_{n-1}\le c_n\,\|g\|_{L^{n,\infty}(\R^n)}
\end{equation*}  
for some dimensional constant $c_n>0$ (for instance, see \cite{S20-New}*{Section 6}), so that one can recover the above inequalities~\eqref{eqi:hardy_bv}, \eqref{eqi:GNS_bv} and~\eqref{eqi:alvino_bv} from~\eqref{eqi:MZ_trace} via known interpolation inequalities in Lorentz spaces.

Motivated by the analogy between $BV$ and $BV^\alpha$ functions, one would be tempted to say that, at least for $n\ge2$, an inequality of the form
\begin{equation}
\label{eqi:buco}
\int_{\R^n}|f|\,d\mu
\le c_{n,\alpha}\,\|\mu\|_{n-\alpha}
\int_{\R^n}|\nabla^\alpha f|\,dx,
\end{equation}
that is, equivalently,
\begin{equation}
\label{eqi:buchetto}
\int_{\R^n} |I_{\alpha} f|\,d\mu
\le
c_{n,\alpha}\,\|\mu\|_{n-\alpha}
\int_{\R^n}|Rf|\,dx,
\end{equation}
may hold for all sufficiently regular functions $f$, see~\cite{S20-New}*{Question~7.1}, where 
\begin{equation*}
R f(x)
=
\pi^{-\frac{n+1}2}\,\Gamma\left(\tfrac{n+1}{2}\right)\,\lim_{\eps\to0^+}\int_{\set*{|y|>\eps}}\frac{y\,f(x+y)}{|y|^{n+1}}\,dy,
\quad
x\in\R^n,
\end{equation*}
is the (vector-valued) \emph{Riesz transform} of~$f$.
Unfortunately, in~\cite{S19}*{Theorem~1.3}, Spector ruled out the validity of~\eqref{eqi:buco}, as well as of~\eqref{eqi:buchetto}, whenever $\alpha\in(0,1)$.

Nonetheless, recalling that $\nabla^\alpha f=\nabla I_{1-\alpha}f$, one may apply the Meyers--Ziemer trace inequality~\eqref{eqi:MZ_trace} to the function $I_{1-\alpha}f$ to get
\begin{equation}
\label{eqi:dolcetto}
\int_{\R^n} |(I_{1-\alpha} f)^\star| \,d\mu
\le
c_n\,\|\mu\|_{n-1}
\int_{\R^n}|\nabla^\alpha f|\,dx.
\end{equation}
Interestingly, inequality~\eqref{eqi:dolcetto} turns out to behave as the mother inequality for the Meyers--Ziemer trace inequality~\eqref{eqi:MZ_trace} as well as for the fractional Hardy inequality~\eqref{eqi:hardy_SS}. 
Indeed, on the one side, taking the limit as $\alpha\to1^-$ in~\eqref{eqi:dolcetto}, then one gets inequality~\eqref{eqi:MZ_trace} back.
On the other side, if one takes $f\ge0$ and $\mu=\frac1{|\,\cdot\,|}\,\Leb{n}$, then one easily recognizes that 
\begin{equation*}
\int_{\R^n}I_{1-\alpha}f\,\frac{dx}{|x|}
=
c_n
\int_{\R^n}I_{1-\alpha}f\,I_{n-1}\,dx
=
c_n
\int_{\R^n}I_{n-\alpha}f\,\,dx
=
c_{n,\alpha}
\int_{\R^n}\frac{f(x)}{|x|^{\alpha}}\,dx,
\end{equation*}
recovering~\eqref{eqi:hardy_SS}.

Having the above observations in mind, our last main result is the following rigorous statement of the inequality~\eqref{eqi:dolcetto}.

\begin{theorem}[Fractional Meyers--Ziemer trace inequality]\label{res:MZ_trace}
Let $\alpha\in(0,1)$ and $p\in\left[1,\frac{n}{1-\alpha}\right)$. There exists a dimensional constant $c_n>0$ such that
\begin{equation}\label{eq:trace}
\int_{\R^n}|(I_{1-\alpha} f)^\star|\,d\mu
\le c_n\,\|\mu\|_{n-1}\,|D^\alpha f|(\R^n)
\end{equation}
for all $f\in BV^{\alpha,p}(\R^n)$ and all $\mu\in\M_{\loc}^+(\R^n)$.
\end{theorem}

As formally observed above, besides providing an alternative route for the proof of \cref{thm:Hardy}, \cref{res:MZ_trace} leads to the following consequences.
Here and in the rest of the paper, we let 
\begin{equation*}
\mathcal H^1(\R^n)
=\set*{f\in L^1(\R^n) : Rf\in L^1(\R^n;\R^n)}
\end{equation*}
be the (real) \emph{Hardy space}, see~\cites{G14-M,S93} for a detailed exposition.

\begin{corollary}[Meyers--Ziemer trace inequalities]
\label{res:MZ_limit}
There exists a dimensional constant $c_n>0$ with the following properties.
\begin{enumerate}[(i)]
\item\label{item:MZ_new} If $f\in BV^{1,p}(\R^n)$ for some $p\in[1,+\infty)$, with $p\le\frac n{n-1}$ if $n\ge2$, and $\mu\in\M_{\loc}^+(\R^n)$, then
\begin{equation}\label{eq:MZ_new}
\int_{\R^n}|f^\star|\,d\mu
\le
c_n\,\|\mu\|_{n-1}\,|Df|(\R^n).
\end{equation} 
\item\label{item:boh} 
If $f\in \mathcal H^1(\R^n)$ and $\mu\in\M_{\loc}^+(\R^n)$, then
\begin{equation}\label{eq:boh}
\int_{\R^n}|(I_1 f)^\star|\,d\mu
\le
c_n\,\|\mu\|_{n-1}\,\|Rf\|_{L^1(\R^n;\,\R^n)}.
\end{equation} 
\end{enumerate}
\end{corollary}

We notice that \cref{res:MZ_limit}\eqref{item:boh} positively answers~\cite{S20-New}*{Question~7.1} in the (solely possible) case $\alpha=1$ and, as well-known, it implies the following stronger version of the Stein--Weiss inequality,
\begin{equation}
\label{eqi:SW}
\|I_1f\|_{L^{\frac n{n-1}}(\R^n)}
\le 
c_n\,\|Rf\|_{L^1(\R^n;\,\R^n)}
\end{equation}
for all $f\in\mathcal H^1(\R^n)$, see~\cite{S20-New}*{Section~1} for a more detailed discussion.
Consequently, once again choosing the measure $\mu=\frac1{|\,\cdot\,|}\,\Leb{n}$,   inequality~\eqref{eq:boh} implies the Hardy-type inequality
\begin{equation}
\label{eqi:hardy_1}
\int_{\R^n}\frac{|I_1f(x)|}{|x|}\,dx
\le 
c_n\,\|Rf\|_{L^1(\R^n;\,\R^n)}
\end{equation} 
whenever $f\in \mathcal H^1(\R^n)$.
Inequality~\eqref{eqi:hardy_1}, in turn, can be also inferred from the Hardy inequality~\eqref{eqi:hardy_bv}, thanks to the continuity of the map $I_1\colon\mathcal H^1(\R^n)\to BV^{1,\frac n{n-1}}(\R^n)$ provided by~\eqref{eqi:SW} (see~\cite{BCCS20}*{Proposition~3.4(i)} for the fractional case $\alpha\in(0,1)$).

\subsection{Organization of the paper}

The paper is organized as follows.
\cref{sec:proof_lim_GG} is dedicated to the proof of  \cref{res:lim_GG_formula}.
In \cref{sec:applications_lim_GG}, we apply it first to prove \cref{thm:IBP_halfspace} and then, in turn, its consequences 
\cref{res:weighted_hardy},
\cref{res:null_function},
\cref{res:chain_sberla_debole}.
Finally, in \cref{sec:MZ_proof}, we prove \cref{res:MZ_trace} and its consequences in \cref{res:MZ_limit}. 

\section{Proof of \texorpdfstring{\cref{res:lim_GG_formula}}{of the limit integration-by-parts formula}}

\label{sec:proof_lim_GG}

In the proof \cref{res:lim_GG_formula}, we take advantage of the following non-local Leibniz rule for $BV^{\alpha,p}$ functions, see~\cite{CS21}*{Theorem~1.1 and Corollary~2.7}. 
For $p \in \left [1, \frac{n}{n-\alpha}\right)$ and $q \in \left ( \frac{n}{\alpha}, + \infty \right ]$ such that  $\frac{1}{p} + \frac{1}{q} = 1$, if $f \in BV^{\alpha,p}(\R^n)$ and $g \in B^{\alpha}_{q, 1}(\R^n)$, then $fg \in BV^{\alpha, r}(\R^n)$ for all $r \in [1, p]$, with $\nabla^{\alpha}_{\rm NL}(f,g) \in L^1(\R^n; \R^n)$ and
\begin{equation}
\label{eq:leibniz_X} 
D^{\alpha}(fg) 
= 
g^\star  D^{\alpha} f 
+ 
f\,
\nabla^{\alpha}g\,
\Leb{n} 
+ 
\nabla^{\alpha}_{\rm NL}(f,g)\,\Leb{n} \quad 
\text{in}\ \M (\R^n; \R^{n}).
\end{equation}
Here and in the rest of the paper, we let
\begin{equation*}
\label{eqi:def_frac_nabla_NL}
\nabla^\alpha_{\NL} (f,g)(x)
=
\mu_{n,\alpha}
\int_{\R^n}\frac{(y-x)(f(y)-f(x))(g(y)-g(x))}{|y-x|^{n+\alpha+1}}\,dy,
\quad
\text{for a.e.}\ x\in\R^n,	
\end{equation*}
be the \emph{non-local fractional $\alpha$-gradient} of the couple $(f,g)$.

\begin{proof}[Proof of \cref{res:lim_GG_formula}] Let $R>0$ be fixed.
Since $\eta_R\in B^\alpha_{q,1}(\R^n)$ with $q\in\left(\frac n\alpha,+\infty\right]$, by the Sobolev Embedding Theorem (see~\cite{AF03}*{Theorem~7.34(c)} and~\cite{L17}*{Theorem~17.52}  for instance) we know that $\eta_R\in C_b(\R^n)$.
Now let $(\rho_\eps)_{\eps>0}$ be a family of standard mollifiers (see~\cite{CS19}*{Section~3.3} for example) and let $g_\eps=\rho_\eps*g$ for all $\eps>0$.
We note that $g_\eps\in\Lip_b(\R^n)$ and $\nabla^\alpha g_\eps=\rho_\eps*\nabla^\alpha g$ for all $\eps>0$, so that 
\begin{equation*}
\lim_{\eps\to0^+}
\int_{\R^n}\eta_R\,f\,\nabla^\alpha g_\eps\,dx
=
\lim_{\eps\to0^+}
\int_{\R^n}\rho_\eps*(f\eta_R)\,\nabla^\alpha g\,dx
=
\int_{\R^n}\eta_R\,f\,\nabla^\alpha g\,dx
\end{equation*}
by the Dominated Convergence Theorem, since 
\begin{equation*}
|\rho_\eps*(f\eta_R)|
\le 
\|f\|_{L^\infty(\R^n)}
\,
\|\eta_R\|_{L^\infty(\R^n)}
\,
\chi_{A_R}
\end{equation*}
for all $\eps>0$ sufficiently small, where $A_R\subset\R^n$ is a bounded set such that $A_R\supset\supp\eta_R$.
Now let $\eps>0$ be fixed. By~\eqref{eq:leibniz_X}, we have that $f\eta_R\in BV^\alpha(\R^n)$, with
\begin{equation*}
D^\alpha(f\eta_R)
=
\eta_R\,D^\alpha f
+
f\,\nabla^\alpha\eta_R\,\Leb{n}
+
\nabla^\alpha_{\NL}(f,\eta_R)\,\Leb{n}
\quad
\text{in}\ \M(\R^n;\R^n).
\end{equation*}
Consequently, by~\cite{CS19-2}*{Proposition~2.7}, we can compute
\begin{align*}
\int_{\R^n}\eta_R\,f\,\nabla^\alpha g_\eps\,dx
&=
-\int_{\R^n}g_\eps\,dD^\alpha (f\eta_R)\\\
&
=
-
\int_{\R^n}\eta_R\, g_\eps\,dD^\alpha f
-
\int_{\R^n}f\,g_\eps\,\nabla^\alpha\eta_R\,dx
-
\int_{\R^n}g_\eps\nabla^\alpha_{\NL}(f,\eta_R)\,dx.
\end{align*}
On the one side, we can estimate
\begin{align*}
\bigg|
\int_{\R^n}f\, g_\eps\nabla^\alpha \eta_R\,dx
\,\bigg|
&\le
\|g_\eps\|_{L^\infty(\R^n)}
\int_{\R^n}|f|\,|\nabla^\alpha\eta_R|\,dx
\\
&\le
\|g\|_{L^\infty(\R^n)}
\,
\|f\|_{L^p(\R^n)}
\,
\|\nabla^\alpha\eta_R\|_{L^q(\R^n;\,\R^n)}
\\
&\le
\mu_{n,\alpha}\,
R^{\frac nq-\alpha}
\,
\|f\|_{L^p(\R^n)}
\,
\|g\|_{L^\infty(\R^n)}
\,
[\eta]_{B^\alpha_{q,1}(\R^n)},
\end{align*}
thanks to~\cite{CS21}*{Corollary~2.3}.
On the other side, in a similar way, we can bound
\begin{align*}
\bigg|
\int_{\R^n}f\,\nabla^\alpha_{\NL}(g_\eps,\eta_R)\,dx
\,\bigg|
&\le
\int_{\R^n}|f|\, \left | \nabla^\alpha_{\NL}(g_\eps,\eta_R) \right | \,dx
\\
&\le
\|f\|_{L^p(\R^n)}
\,
\|\nabla^\alpha_{\NL}(g_\eps,\eta_R)\|_{L^q(\R^n;\,\R^n)}
\\
&\le
2\mu_{n,\alpha}
\,
\|f\|_{L^p(\R^n)}
\,
\|g_\eps\|_{L^\infty(\R^n)}
\,
[\eta_R]_{B^\alpha_{q,1}(\R^n)}
\\
&\le
2\mu_{n,\alpha}
\,
R^{\frac nq-\alpha}
\,
\|f\|_{L^p(\R^n)}
\,
\|g\|_{L^\infty(\R^n)}
\,
[\eta]_{B^\alpha_{q,1}(\R^n)}
\end{align*}
thanks to~\cite{CS21}*{Corollary~2.7}.
Therefore, thanks to these estimates (which are uniform in $\eps$), we get the limit
\begin{equation*}
\lim_{R\to+\infty}
\sup_{\eps>0}
\bigg|
\int_{\R^n}f\, g_\eps\nabla^\alpha \eta_R\,dx
\,\bigg|
+
\bigg|
\int_{\R^n}f\,\nabla^\alpha_{\NL}(g_\eps,\eta_R)\,dx
\,\bigg|
=0.
\end{equation*}
Now we need to show that
\begin{equation}
\label{eq:fulmine}
\lim_{\eps\to0^+}
\int_{\R^n}\eta_R\, g_\eps\,dD^\alpha f
=
\int_{\R^n}\eta_R\, g^\star \,dD^\alpha f.
\end{equation}
Indeed, since $f\in BV^{\alpha,\infty}(\R^n)$, by~\cite{CSS21}*{Theorem~1} we have that $|D^\alpha f|\ll\Haus{n-\alpha}$. 
Moreover, being $g\in W^{\alpha,1}_{\loc}(\R^n)$, by~\cite{PS20}*{Proposition~3.1} we can infer that  
\begin{equation*}
\lim_{\eps\to0^+}
g_\eps(x)
=
\lim_{\eps\to0^+}
\rho_\eps*g(x)
=  \lim_{r\to0^+} \aint_{B_r(x)} g(y)\,dy =
g^\star(x)
\quad
\text{for}\ \Haus{n-\alpha}\text{-a.e.}\ x\in\R^n,
\end{equation*}
so that~\eqref{eq:fulmine} immediately follows by the Dominated Convergence Theorem (with respect to the finite measure $|D^\alpha f|$).
Finally, since 
\begin{equation*}
\lim_{R\to+\infty}
\eta_R(x)
=
\lim_{R\to+\infty}
\eta\left(\tfrac xR\right)
=
\eta(0)
\quad
\text{for all}\ x\in\R^n,
\end{equation*}
by the Dominated Convergence Theorem (with respect to the finite measure $|g^\star|\,|D^\alpha f|$) we conclude that
\begin{equation*}
\lim_{R\to+\infty}
\int_{\R^n}\eta_R\, g^\star \,dD^\alpha f
=
\eta(0)
\int_{\R^n} g^\star \,dD^\alpha f
\end{equation*}
and the proof is complete.
\end{proof}

\section{Hardy inequalities and failure of the chain rule}

\label{sec:applications_lim_GG}

\subsection{Integration by parts on half-spaces}

We begin with the proof of the formula for the fractional gradient of the characteristic function of a half-space.

\begin{proof}[Proof of \cref{prop:halfspace_nabla_alpha}]
By the translation invariance of the fractional gradient (recall~\cite{Sil19}*{Theorem~2.2}), we have
\begin{equation*}
\nabla^{\alpha} \chi_{H_{\nu}^{+}(x_0)}(x) = \nabla^{\alpha} \chi_{H_{\nu}^{+}}(x - x_0)
\end{equation*}
for all $x\in\R^n$ and so we can assume $x_0 = 0$ without loss of generality.
Since $\chi_{H_\nu^+}\in BV_{\loc}(\R^n)\cap L^\infty(\R^n)$ and clearly 
\begin{equation*}
|D \chi_{H^{+}_{\nu}}|(\partial B_{R}) = \Haus{n - 1}(H_{\nu} \cap \partial B_{R}) = 0
\end{equation*}
for all $R > 0$, by~\cite{CS19-2}*{Proposition~3.5} we get $\nabla^\alpha\chi_{H_\nu^+}  \in L^1_{\loc}(\R^n;\R^n)$ and
\begin{align} 
\int_{\R^{n}} \varphi \cdot \nabla^{\alpha} \chi_{H^{+}_{\nu}} \, dx & = \lim_{R \to + \infty} \int_{\R^{n}} \varphi \cdot I_{1 - \alpha} ( \chi_{B_{R}} D \chi_{H^{+}_{\nu}}) \, dx \nonumber \\
& = \lim_{R \to + \infty} \nu\cdot\int_{\R^{n}} \varphi\, I_{1 - \alpha} ( \chi_{B_{R}} \Haus{n - 1} \res H_{\nu}) \, dx \label{eq:dual_repr_halfspace_lim}
\end{align}
for all $\varphi \in \Lip_c(\R^n; \R^n)$.
By the Monotone Convergence Theorem, we get
\begin{equation*}
\lim_{R\to+\infty}
I_{1 - \alpha} ( \chi_{B_{R}} \Haus{n - 1} \res H_{\nu})(x) 
=
I_{1 - \alpha} (\Haus{n - 1} \res H_{\nu})(x)
\end{equation*} 
for $\Leb{n}$-a.e.\ $x \in \R^{n}$.
We now claim that  
\begin{equation}
\label{eq:Riesz_pot_halfspace_alpha}
I_{1 - \alpha}(\Haus{n - 1} \res H_{\nu})(x)
=
\frac{\mu_{1,\alpha}}{\alpha}
\,
\frac{1}{|x\cdot \nu|^\alpha} \ \text{ for all } x \notin H_{\nu},
\end{equation}
which defines a function in $L^{1}_{\rm loc}(\R^{n})$.
The case $n = 1$ is easy.
For $n \ge 2$, we argue as follows.
Let $\mathcal{R} \in\SO(n)$ be such that $\mathcal{R} \nu = \e_{1}$, so that 
$(\mathcal{R}x)_{1} = (\mathcal{R}x) \cdot \e_{1} = x \cdot \transp{\mathcal{R}} \e_{1} = x \cdot \nu$. 
By simple changes of variables, we get
\begin{align*}
\int_{H_{\nu}} \frac{d \Haus{n - 1}(y)}{|y - x|^{n + \alpha - 1}} 
&=  
\int_{H_{\e_{1}}} \frac{d \Haus{n - 1}(y)}{|y - \mathcal R x|^{n + \alpha - 1}} 
\\
& = 
\int_{\R^{n - 1}} \frac{dy_2\,\cdots\, dy_n}{\big( ({\mathcal R} x)_{1}^{2} + \sum_{j = 2}^{n} (y_{j} - ({\mathcal R} x)_{j})^{2} \big)^{\frac{n + \alpha - 1}{2}}} 
\\
& = 
\int_{\R^{n - 1}}
\frac{1}{|({\mathcal R}x)_{1}|^{\alpha}} 
\, \frac{dy_2\,\cdots\, dy_n}{\big(1 + |(y_2,\dots,y_n)|^{2}\big)^{\frac{n + \alpha - 1}{2}}} 
\\
& = 
\frac{(n - 1)\, \omega_{n - 1}}{|x \cdot \nu|^{\alpha}}  \int_{0}^{+\infty} \frac{\rho^{n - 2}}{(1 + \rho^{2})^{\frac{n + \alpha - 1}{2}}} \, d \rho
\end{align*}
whenever $x \notin H_{\nu}$.
By known properties of the Gamma function, it is not difficult to recognize that 
\begin{align*}
\int_{0}^{+\infty} \frac{\rho^{n - 2}}{(1 + \rho^{2})^{\frac{n + \alpha - 1}{2}}} \, d \rho 
= \frac{\Gamma\left(\frac{\alpha}{2}\right) \Gamma\left(\frac{n - 1}{2}\right)}{2\, \Gamma\left(\frac{n + \alpha - 1}{2}\right)},
\end{align*}
so that 
\begin{align*}
I_{1 - \alpha}(\Haus{n - 1} \res H_{\nu})(x) 
& = 
\frac{\mu_{n,\alpha}}{(n + \alpha - 1)} 
\frac{\Gamma\left(\frac{\alpha}{2}\right) \Gamma\left(\frac{n - 1}{2}\right)}{2\, \Gamma\left(\frac{n + \alpha - 1}{2}\right)} \frac{(n - 1)\, \omega_{n - 1}}{|x \cdot \nu|^{\alpha}} \nonumber 
= 
\frac{\mu_{1,\alpha}}{\alpha}
\frac{1}{|x \cdot \nu|^{\alpha}} 
\end{align*}
whenever $x\notin H_\nu$, proving~\eqref{eq:Riesz_pot_halfspace_alpha}. 
Therefore, we can apply the Dominated Convergence Theorem in~\eqref{eq:dual_repr_halfspace_lim} to obtain
\begin{equation*}
\int_{\R^{n}} \varphi \cdot \nabla^{\alpha} \chi_{H^{+}_{\nu}} \, dx 
=
\nu\cdot\int_{\R^{n}} \varphi \,I_{1 - \alpha}(\Haus{n - 1} \res H_{\nu}) \, dx
=
\frac{\mu_{1,\alpha}}{\alpha}
\,
\nu\cdot
\int_{\R^n}\frac{\phi(x)}{|x\cdot \nu|^\alpha}\,dx
\end{equation*}
and the conclusion immediately follows.
\end{proof}

Having \cref{prop:halfspace_nabla_alpha} at disposal, we can easily deduce the limit Gauss--Green formula on half-spaces. 

\begin{proof}[Proof of \cref{thm:IBP_halfspace}]
The validity of~\eqref{eq:IBP_halfspace} is an immediate consequence of \cref{res:lim_GG_formula}, since $\chi_{H^+_\nu(x_0)}\in BV_{\loc}(\R^n)\cap L^\infty(\R^n)$ with $\nabla^\alpha\chi_{H^+_\nu(x_0)}\in L^1_{\loc}(\R^n;\R^n)$ thanks to \cref{prop:halfspace_nabla_alpha}.
For the proof of~\eqref{eq:IBP_halfspace_special}, we can simply choose $\eta\in C^\infty_c(\R^n)$ such that $0\le\eta\le1$ and $\eta(x)=1$ for $x\in B_1$, so that, arguing component-wise, 
\begin{equation*}
\lim_{R\to+\infty}
\int_{\R^n}\eta_R(x)\,\frac{f(x)\,\nu}{|(x-x_0)\cdot\nu|^\alpha}\,dx
=
\int_{\R^n}\frac{f(x)\,\nu}{|(x-x_0)\cdot\nu|^\alpha}\,dx
\end{equation*} 
either trivially if $\supp f$ is bounded, or by the Monotone Convergence Theorem if $f$ has constant sign. Thus, the proof is complete.
\end{proof}

\subsection{Fractional Hardy inequalities}

We can now deal with the proofs of the fractional Hardy inequalities in \cref{res:hardy_mezza}, \cref{res:weighted_hardy} and \cref{thm:Hardy}.

\begin{proof}[Proof of \cref{res:hardy_mezza}]
Let $(\rho_\eps)_{\eps>0}$ be a family of standard mollifiers and set $f_\eps=\rho_\eps*f$ for all $\eps>0$.
Clearly, $f_\eps\in BV^{\alpha,p}_+(\R^n)\cap L^\infty(\R^n)$, so that \eqref{eq:IBP_halfspace_special} implies 
\begin{equation*}
\frac{\mu_{1,\alpha}}{\alpha}
\int_{\R^n} \frac{f_\eps (x)}{|(x - x_0)\cdot\nu|^{\alpha}} \, dx \le |D^\alpha f_\eps|\big(H_\nu^+(x_0)\big)
\end{equation*}
for all $x_0 \in \R^n$ and $\nu\in\Sph^{n-1}$. 
On the left-hand side, we employ Fatou's Lemma to obtain
\begin{equation*}
\liminf_{\eps \to 0^+ } \frac{\mu_{1,\alpha}}{\alpha}
\int_{\R^n} \frac{f_\eps (x)}{|(x - x_0)\cdot\nu|^{\alpha}} \, dx \ge \frac{\mu_{1,\alpha}}{\alpha}
\int_{\R^n} \frac{f(x)}{|(x - x_0)\cdot\nu|^{\alpha}} \, dx.
\end{equation*} 
As for the right-hand side, thanks to \cite{AFP00}*{Theorem~2.2(b)} we notice that 
$$\limsup_{\eps \to 0^+} |D^\alpha f_\eps|\big(H_\nu^+(x_0)\big) \le \limsup_{\eps \to 0^+} |D^\alpha f|\big(H_\nu^+(x_0) + B_\eps \big) \le |D^\alpha f|\big(\closure[0]{H_\nu^+(x_0)}\big)
$$
and this proves \eqref{eq:hardy_mezza}.
\end{proof}

\begin{proof}[Proof of \cref{res:weighted_hardy}]
At first, let us also assume that $f\in BV^{\alpha, p}_+(\R^n) \cap L^\infty(\R^n)$.
Let $r>0$ be fixed. 
Choosing $x_0+r\nu$ in place of $x_0$ in~\eqref{eq:IBP_halfspace_special} and taking the integral average on $\Sph^{n-1}$, we get
\begin{align*}
\frac{\mu_{1,\alpha}}{\alpha}
\,
\aint_{\Sph^{n-1}}
\int_{\R^n} \frac{f(x)}{|(x - x_0)\cdot\nu-r|^{\alpha}} \, dx
\,d\Haus{n-1}(\nu)
&=
- 
\aint_{\Sph^{n-1}}
\nu\cdot D^{\alpha} f(H_{\nu}^+(x_0+r\nu))
\,d\Haus{n-1}(\nu)
\\
&\le 
|D^\alpha f|\big(\R^n\setminus B_r(x_0)\big).
\end{align*}
By Tonelli's Theorem, we can compute
\begin{align*}
\aint_{\Sph^{n-1}}
\int_{\R^n} \frac{f(x)}{|(x - x_0)\cdot\nu-r|^{\alpha}} \, dx
\,d\Haus{n-1}(\nu)
&=
\int_{\R^n}
f(x)\,
\aint_{\Sph^{n-1}}
\frac{d\Haus{n-1}(\nu)}{|(x - x_0)\cdot\nu-r|^{\alpha}} \,dx
\\
&=
\int_{\R^n}
f(x)\,w_{n,\alpha}(|x-x_0|,r)\,dx,
\end{align*}
where in the last inequality we exploited the formula proved in~\cite{G14-C}*{Section~D.3} for $n\ge2$ (the case $n=1$ being trivial). 

Now let $f\in BV^{\alpha,p}_+(\R^n)$ be possibly unbounded.
Let $(\rho_\eps)_{\eps>0}$ be a family of standard mollifiers and set $f_\eps=\rho_\eps*f$ for all $\eps>0$.
Clearly, $f_\eps\in BV^{\alpha,p}_+(\R^n)\cap L^\infty(\R^n)$, so that
\begin{equation*}
\int_{\R^n}
f_\eps(x)\,w_{n,\alpha}(|x-x_0|,r)\,dx
\le 
|D^\alpha f_\eps|\big(\R^n\setminus B_r(x_0)\big)
\end{equation*}
for all $\eps>0$. 
On the one side, we have
\begin{equation*}
\liminf_{\eps\to0^+}
\int_{\R^n}
f_\eps(x)\,w_{n,\alpha}(|x-x_0|,r)\,dx
\ge 
\int_{\R^n}
f(x)\,w_{n,\alpha}(|x-x_0|,r)\,dx
\end{equation*}
by Fatou's Lemma. 
On the other side, thanks to~\cite{CSS21}*{Theorem~4} and~\cite{AFP00}*{Theorem~2.2(b)}, we can estimate
\begin{equation*}
|D^\alpha f_\eps|\big(\R^n\setminus {B_r(x_0)}\big)
\le
(\rho_\eps*|D^\alpha f|)\big(\R^n\setminus {B_r(x_0)}\big)
\le
|D^\alpha f|\big(\R^n\setminus B_{r-\eps}(x_0)\big)
\end{equation*} 
for all $\eps\in(0,r)$.
Consequently,
we get
\begin{align*}
\int_{\R^n}
f(x)\,w_{n,\alpha}(|x-x_0|,r)\,dx
\le
\lim_{\eps\to0^+}
|D^\alpha f|\big(\R^n\setminus B_{r-\eps}(x_0)\big)
=
|D^\alpha f|\big(\R^n\setminus B_{r}(x_0)\big)
\end{align*}
by monotonicity and the proof is complete.
\end{proof}

\begin{proof}[Proof of \cref{thm:Hardy}]
At first, let $n \ge 2$ and $f\in BV^{\alpha,p}(\R^n)$ with $p\in\left[1,\frac n{1-\alpha}\right)$. Up to a translation, we can assume $x_0 = 0$.
Let $(\rho_\eps)_{\eps>0}$ be a family of standard mollifiers and let $f_\eps=\rho_\eps*f$ for all $\eps>0$.
By~\cite{CSS21}*{Theorem~4}, we know that $f_\eps\in BV^{\alpha,p}(\R^n)$ and $\nabla^\alpha f_\eps=\rho_\eps*D^\alpha f$ for all $\eps>0$.
Moreover, thanks to~\cite{CS19}*{Lemma~3.28(i)} and~\cite{CSS21}*{Proposition~4(i)}, we have $\nabla^\alpha f_\eps=\nabla I_{1-\alpha}f_\eps\in L^1(\R^n;\R^n)$.
Therefore the conclusion follows by applying~\eqref{eqi:hardy_spector} to $f_\eps$ and then passing to the limit as $\eps\to0^+$ via Fatou's Lemma and \cite{CSS21}*{Theorem~4}. Let now $n = 1$, $f \in BV^{\alpha, p}_{+}(\R)$ with $p \in \left [1, \frac{1}{1 - \alpha} \right )$, and $f_\eps=\rho_\eps*f$ for all $\eps>0$. Clearly, $f_\eps\ge0$, so that we may employ \cref{res:hardy_mezza} to get
\begin{equation} \label{eq:Hardy_half_line}
\frac{\mu_{1,\alpha}}{\alpha}
\int_{\R} \frac{f_\eps(x)}{|x - x_0|^{\alpha}} \, dx = \frac{\mu_{1,\alpha}}{\alpha}
\int_{\R} \frac{f_\eps(x)}{|(x - x_0)\cdot\nu|^{\alpha}} \, dx \le |D^\alpha f_{\eps}|(H_\nu^+(x_0))
\end{equation}
for all $x_0 \in \R$ and $\nu\in \{\pm 1\}$, since $f_\eps \in BV^{\alpha, \infty}(\R)$, and so $|D^\alpha f_\eps| \ll \Haus{1 - \alpha}$ by~\cite{CSS21}*{Theorem~1}. Hence, if we substitute $\nu$ with $-\nu$ in \eqref{eq:Hardy_half_line} and then add the two inequalities, we get
\begin{equation*}
\frac{2 \mu_{1,\alpha}}{\alpha}
\int_{\R} \frac{f_\eps(x)}{|x - x_0|^{\alpha}} \, dx \le |D^\alpha f_{\eps}|(\R) \le |D^\alpha f|(\R).
\end{equation*}
Thus, we can pass to the limit as $\eps \to 0^+$ exploiting again Fatou's Lemma. In order to prove the optimality of the constant $c_{1, \alpha} = \frac{2 \mu_{1,\alpha}}{\alpha}$, we choose $f = \chi_{(x_0 -1, x_0 + 1)}$, so that
\begin{equation*}
\int_{\R} \frac{\chi_{(x_0 -1, x_0 + 1)}(x)}{|x - x_0|^{\alpha}} \, dx = \frac{2}{1 - \alpha}.
\end{equation*}
Since 
\begin{equation*}
|D^\alpha \chi_{(x_0 -1, x_0 + 1)}|(\R) = \frac{4 \mu_{1, \alpha}}{\alpha (1 - \alpha)}
\end{equation*}
thanks to \cite{CS19}*{Example 4.11}, we get the optimality of $c_{1, \alpha}$ and the proof is complete.
\end{proof}

\begin{remark}
Let $\alpha \in (0, 1)$ and $p \in \left [1, \frac{n}{n - \alpha} \right )$.
Arguing as in the second part of the proof of \cref{thm:Hardy}, it is possible to show that
\begin{equation*}
\frac{2 \mu_{1,\alpha}}{\alpha}
\int_{\R^n} \frac{f(x)}{|x - x_0|^{\alpha}} \, dx \le |D^\alpha f|(\R^n) 
\end{equation*}
for all $f \in BV^{\alpha, p}_{+}(\R^n)$. Combining this with \cite{SS18}*{Theorem 1.2}, we deduce that 
\begin{equation*}
c_{n, \alpha}
\int_{\R^n} \frac{f(x)}{|x - x_0|^{\alpha}} \, dx \le |D^\alpha f|(\R^n) \text{ for all } f \in BV^{\alpha, p}_{+}(\R^n),
\end{equation*}
where
\begin{equation*}
c_{n, \alpha} = \max \left \{\frac{2 \mu_{1,\alpha}}{\alpha}, \gamma_{n, \alpha} \right \}
\end{equation*}
and
\begin{equation*}
\gamma_{n, \alpha} = \frac{2^{\alpha} \Gamma \left ( \frac{\alpha}{2} \right ) \Gamma \left ( \frac{n+1}{2} \right ) }{\pi^{1 - \frac{\alpha}{2}} \Gamma \left ( \frac{n - \alpha}{2} \right )}.
\end{equation*}
However, for $n \ge 3$, one can see that $\gamma_{n, \alpha} > \frac{2 \mu_{1,\alpha}}{\alpha}$ for all $\alpha \in (0, 1)$.
To the best of our knowledge, it is not known whether $c_{n, \alpha}$ is the optimal constant for some $n \ge 2$ and $\alpha \in (0, 1)$.
\end{remark}

\subsection{Failure of the fractional chain rule}

We begin with the proof of the rigidity property contained in \cref{res:null_function}. 

\begin{proof}[Proof of \cref{res:null_function}]
If $\supp|D^\alpha f|$ is bounded, then $|D^\alpha f|\big(\R^n\setminus B_r\big)=0$ for some $r>0$. 
Hence, by \cref{res:weighted_hardy}, we must have $f=0$ $\Leb{n}$-a.e.\ in $\R^n$, being \mbox{$w_{n,\alpha}>0$}. 
If, instead, $|D^\alpha f|\big(\closure[0]{H_\nu^+(x_0)}\big) = 0$ or $f\in L^\infty(\R^n)$ and $D^\alpha f(H^+_\nu(x_0))=0$ for some $x_0\in\R^n$ and $\nu\in\Sph^{n-1}$, then we similarly conclude by~\eqref{eq:IBP_halfspace_special} in \cref{thm:IBP_halfspace} and \cref{res:hardy_mezza}.
\end{proof}

We can now end this section by showing the failure of the fractional chain rule. 
Here and in the following, we let
\begin{equation*}
(- \Delta)^{\frac{\beta}{2}} f(x) 
= 
\nu_{n, \beta}
\int_{\R^n} \frac{f(x + y) - f(x)}{|y|^{n + \beta}}\,dy,
\quad
\text{for a.e.}\ x\in\R^n,
\end{equation*}
be the \emph{fractional Laplacian} of order $\beta \in(0,1)$ of the function $f \in W^{\beta, 1}(\R^n)$,
where
\begin{equation*}
\nu_{n,\beta}=2^\beta \pi^{-\frac n2}\frac{\Gamma\left(\frac{n+\beta}{2}\right)}{\Gamma\left(-\frac{\beta}{2}\right)}.
\end{equation*}
Note that $(- \Delta)^{\frac{\beta}{2}}\colon W^{\beta, 1}(\R^n)\to L^1(\R^n)$ is continuous (see~\cite{CS19}*{Section~3.10} for a more detailed discussion). 
In particular, this operator is well-posed on $BV$ functions.

\begin{proof}[Proof of \cref{res:chain_sberla_forte}]
Let $Q_1=(-1,1)^n$. 
We consider the function $f=(-\Delta)^{\frac{1-\alpha}2}\chi_{Q_1}$, that is,
\begin{equation}
\label{eq:magic_f}
f(x)
=
\nu_{n,1-\alpha}
\left(
-\chi_{Q_1}(x)
\int_{\R^n\setminus Q_1}
\frac{1}{|y-x|^{n+1-\alpha}}\,dy
+
\chi_{\R^n\setminus Q_1}(x)
\int_{Q_1}
\frac{1}{|y-x|^{n+1-\alpha}}\,dy
\right),
\end{equation} 
for $x\in\R^n\setminus\de Q_1$.
Thanks to~\cite{CS19}*{Lemma~3.28(ii)}, we know that  $f\in BV^\alpha(\R^n)$ with $D^\alpha f=D\chi_{Q_1}$.
By~\cite{CSS21}*{Theorem~6}, we also have that $f\in BV^{\alpha,p}(\R^n)$ for all $p\in\left[1,\frac n{n-\alpha}\right)$. 
Now let $\Phi\in\Lip_b(\R)$ be such that $\Phi(0)=0$ and $\Phi\ge0$ and assume that $\Phi(f)\in BV^\alpha(\R^n)$ with $$\supp |D^\alpha\Phi(f)| \subset \supp |D^\alpha f|= \supp |D\chi_{Q_1}| = \partial Q_1.$$
Note that, again by~\cite{CSS21}*{Theorem~6}, $\Phi(f)\in BV^{\alpha,p}_+(\R^n)$ for all $p\in\left[1,\frac n{n-\alpha}\right)$. 
Consequently, $\supp|D^\alpha\Phi(f)|$ is compact, so that $\Phi(f)\equiv0$ thanks to \cref{res:null_function}. 
Since $\nu_{n,1-\alpha} < 0$, we observe that $f(x)\to0^-$ as $|x|\to+\infty$ and, moreover, 
\begin{align*}
\liminf_{t\to1^+}
\int_{Q_1}
\frac{1}{|y-t\e_1|^{n+1-\alpha}}\,dy
&\ge
\int_{Q_1}\frac{dy}{|y-\e_1|^{n+\alpha-1}}
\\
&\ge
\sup_{\eps\in(0,1)}
\int_{(-1,1-\eps)\times(-\eps,\eps)^{n-1}}\frac{dy}{|y-\e_1|^{n+\alpha-1}}
\\
&\ge
c_{n,\alpha}
\sup_{\eps\in(0,1)}
\int_{(-1,1-\eps)\times(-\eps,\eps)^{n-1}}\frac{dy}{|y_1-1|^{n+\alpha-1}}
\\
&=
c_{n,\alpha}
\sup_{\eps\in(0,1)}
\eps^{n-1}\left(\eps^{\alpha-n}-2^{\alpha-n}\right)
=+\infty,
\end{align*}
thanks to Fatou's Lemma.
As a consequence, $f(\R^n)\supset(-\infty, 0)$ and thus $\Phi(t)=0$ for all $t\in (-\infty, 0)$.
Replacing $f$ with $-f$, we also get that $\Phi(t)=0$ for all $t\in(0, +\infty)$, proving that $\Phi\equiv0$ and the proof is complete.
\end{proof}

\begin{proof}[Proof of \cref{res:chain_sberla_debole}]
By~\cite{CS19}*{Theorem~3.26}, we know that $f_\alpha\in BV^\alpha(\R)$. 
We claim that $|f_\alpha| \notin BV^\alpha(\R)$.
By contradiction, if $|f_\alpha| \in BV^\alpha(\R)$, then~\cref{thm:Hardy} implies that
\begin{equation}
\label{eq:furetto} 
c_\alpha\int_{\R} \frac{|f_{\alpha}(x)|}{|x|^{\alpha}} \, dx 
\le 
|D^\alpha |f_{\alpha}||(\R) < + \infty.
\end{equation}
However, for $x \in (0, 1)$, we have
\begin{equation*}
\frac{|f_{\alpha}(x)|}{|x|^{\alpha}} = |\mu_{1, - \alpha}| \left ( \frac{1}{x} + \frac{(1- x)^{\alpha - 1}}{x^{\alpha}} \right),
\end{equation*}
contradicting~\eqref{eq:furetto} and the proof is complete.
\end{proof}

\section{Fractional Meyers--Ziemer trace inequality}

\label{sec:MZ_proof}

We begin by noticing that, somehow formulating in a more rigorous way the ideas sketched in the introduction, one can prove \cref{res:MZ_trace} by directly applying the standard Meyers--Ziemer trace inequality~\eqref{eqi:MZ_trace} to the function $u=I_{1-\alpha} f$ whenever $f\in BV^{\alpha,p}(\R^n)$ with $p\in\left(1,\frac{n}{1-\alpha}\right)$, since
\begin{equation*}
I_{1-\alpha}\colon BV^{\alpha,p}(\R^n)\to BV^{1,\frac n{n-1}}(\R^n)
\end{equation*}
with $Du=D^\alpha f$ in~$\M(\R^n;\R^n)$,
thanks to~\cite{CSS21}*{Proposition~4(i)}.
In the case $p=1$, we only have 
\begin{equation*}
I_{1-\alpha}\colon BV^\alpha(\R^n)\to BV^{1,q}(\R^n)
\end{equation*}
for all $q\in\left(\frac{n}{n-1+\alpha},\frac n{n-1}\right)$ with $Du=D^\alpha f$ in~$\M(\R^n;\R^n)$, in virtue of~\cite{CS19}*{Remark~3.29} and~\cite{CSS21}*{Theorem~6}, but this is still enough in order to directly exploit~\eqref{eqi:MZ_trace}. 

Below, we instead outline a direct argument showing that the very same line of reasoning used in~\cite{MZ77} (see~\cite{S20-New}*{Section~7} for a more detailed explanation) to prove~\eqref{eqi:MZ_trace} works as well for proving \cref{res:MZ_trace}.

\begin{proof}[Proof of \cref{res:MZ_trace}]
Let $f\in BV^{\alpha,p}(\R^n)$ for some $p\in\left[1,\frac n{1-\alpha}\right)$.
Let $(\rho_\eps)_{\eps>0}$ be a family of standard mollifiers and let $f_\eps=\rho_\eps*f$ for all $\eps>0$.
By~\cite{CSS21}*{Theorem~4}, we know that $f_\eps\in BV^{\alpha,p}(\R^n)$ with $\nabla^\alpha f_\eps=\rho_\eps*D^\alpha f$ for all $\eps>0$.
Now let $u_\eps=I_{1-\alpha} f_\eps$ for all $\eps>0$.
By what we have just observed above, it is not difficult to see that $u_\eps\in BV^{1,q}(\R^n)\cap C^\infty(\R^n)$ for some $q\in\left(\frac n{n-1+\alpha},\frac n{n-1}\right]$ with
\begin{equation*}
|\nabla u_\eps|
=
|\nabla^\alpha f_\eps|
\le 
\rho_\eps*|D^\alpha f|
\quad 
\text{in}\ L^1(\R^n).
\end{equation*}
Therefore, we can estimate
\begin{equation*}
\int_{\R^n}|\nabla u_\eps|\,dx
=
\int_{\R^n}|\nabla^\alpha f_\eps|\,dx
\le
|D^\alpha f|(\R^n)
<+\infty
\end{equation*} 
and, moreover, the set
\begin{equation*}
E_t^\eps=\set*{x\in\R^n : |u_\eps(x)|>t}
\end{equation*}
is open with finite perimeter for a.e.~$t>0$. Since
\begin{align*}
\frac{|E^\eps_t\cap B_r(x)|}{|B_r(x)|}
\le\frac{\min\set*{|E^\eps_t|,|B_r(x)|}}{|B_r(x)|}
\end{align*}
and
\begin{align*}
|E^\eps_t|=|\set*{x\in\R^n : |I_{1-\alpha} f_\eps(x)|>t}|
\le c_{n,\alpha,p}\left(\frac{\|f\|_{L^p(\R^n)}}{t}\right)^{\frac{np}{n-(1-\alpha)p}}<+\infty
\end{align*}
by the Hardy--Littlewood--Sobolev inequality (see~\cite{G14-M}*{Theorem~1.2.3} for instance), for each given $x\in E^\eps_t$ the function
\begin{equation*}
r\mapsto\frac{|E_t^\eps\cap B_r(x)|}{|B_r(x)|}
\end{equation*}
is continuous, equals~$1$ for small~$r>0$ (since~$E_t^\eps$ is open) and tends to zero as $r\to+\infty$. 
Thus, reasoning exactly as in~\cite{S19}*{Section~6}, via a routine Vitali covering argument we can estimate
\begin{align*}
\mu(E_t^\eps)\le c_n\,\|\mu\|_{n-1}\, |D\chi_{E_t^\eps}|(\R^n)
\end{align*}
for a.e.~$t>0$, where $c_n>0$ is a dimensional constant. 
Therefore, by the coarea formula and the chain rule for functions with bounded variation, we can estimate
\begin{align*}
\int_{\R^n}|u_\eps|\,d\mu
&=
\int_\R\mu(E_t^\eps)\,dt
\\
&\le c_n\,\|\mu\|_{n-1}\int_\R|D\chi_{E_t^\eps}|(\R^n)\,dt
\\
&\le
c_n\,\|\mu\|_{n-1}\int_{\R^n}|\nabla u_\eps|\,dx
\\
&\le
c_n\,\|\mu\|_{n-1}\,|D^\alpha f|(\R^n)
\end{align*}
for all $\eps>0$.
Now, assuming $\|\mu\|_{n-1}<+\infty$ without loss of generality, it is standard to see that $\mu\ll\Haus{n-1}$, see~\cite{PS20} and the references therein for a more detailed discussion.
Therefore, since
\begin{equation*}
\lim_{\eps\to0^+}
u_\eps(x)
=
\lim_{\eps\to0^+}
\rho_\eps*u(x)
=
u^\star(x)
\quad
\text{for}\ \Haus{n-1}\text{-a.e.}\ x\in\R^n
\end{equation*} 
(see~\cite{EG15}*{Section~5.9} for instance),
by the Fatou's Lemma we conclude that 
\begin{equation*}
\int_{\R^n}|(I_{1-\alpha}f)^\star|\,d\mu
\le 
\liminf_{\eps\to0^+}
\int_{\R^n}|u_\eps|\,d\mu
\le
c_n\,\|\mu\|_{n-1}\,|D^\alpha f|(\R^n)
\end{equation*} 
and the proof is complete.
\end{proof}

We now conclude our paper with the proof of \cref{res:MZ_limit}.

\begin{proof}[Proof of \cref{res:MZ_limit}]
The validity of~\eqref{item:MZ_new} for any $f\in C^\infty_c(\R^n)$ follows directly from \cref{res:MZ_trace} combined with the asymptotic analysis obtained in~\cite{CS19-2}.
For a general $f\in BV^{1,\frac{n}{n-1}}(\R^n)$, one just needs to perform a routine approximation argument thanks to~\cite{CSS21}*{Proposition~1}.
The validity of~\eqref{item:boh} follows in a similar way, this time relying on the asymptotic analysis carried out in~\cite{BCCS20}.
\end{proof}


\end{document}